\documentclass[10pt,twocolumn,twoside]{IEEEtran}
\usepackage{amsmath}
\usepackage{amssymb}
\usepackage{amsfonts}
\usepackage{bbm}
\usepackage{color}
\usepackage{enumerate}	
\usepackage{wrapfig}

\usepackage{amssymb,amsfonts,amsmath,mathrsfs,mathtools,mathabx}
\usepackage{graphicx,epsfig}
\usepackage{bm}
\usepackage[margin=1in,papersize={8.5in,11in}]{geometry}
\usepackage{times}
\usepackage{multirow}
\usepackage{color}
\usepackage{xcolor}
\newtheorem{theorem}{Theorem}
\newtheorem{remark}{Remark} 
\newtheorem{assumption}{Assumption}

\newtheorem{lemma}{Lemma}

\newcommand\numberthis{\addtocounter{equation}{1}\tag{\theequation}}
\begin{document}

\title{Defense Against Adversarial Swarms with Parameter Uncertainty}
%
%
%

\author{Claire~Walton,
        Isaac~Kaminer,
        Qi~Gong,
        Abram.~H.~Clark
        and~Theodoros~Tsatsanifos
\thanks{Claire Walton is with the Departments of Mathematics and Electrical \& Computer Engineering, University of Texas at San Antonio, San Antonio, TX 78249 USA (e-mail:~claire.walton@utsa.edu).}
\thanks{Isaac Kaminer and Theodoros Tsatsanifos are are with the Department of Mechanical and Aerospace Engineering,
		Naval Postgraduate School, Monterey, CA, 93943 USA (e-mail:~kaminer@nps.edu;~theodoros.tsatsanifos.gr@nps.edu).}
\thanks{Qi Gong is with the Department of Applied Mathematics, University of California Santa Cruz, Santa Cruz, CA 95064 USA (e-mail:~qigong@ucsc.edu).}
\thanks{Abe Clark is with the Department of Physics, Naval Postgraduate School (e-mail:~abe.clark@nps.edu).}
\thanks{Manuscript received on February 28, 2020.}}

\maketitle

\begin{abstract}
This paper addresses the problem of optimal defense of a High Value Unit against a large-scale swarm attack. We show that the problem can be cast in the framework of uncertain parameter optimal control and derive a consistency result for the dual problem of this framework. We show that the dual can be computed numerically and apply these numerical results to derive optimal  defender strategies against a 100 agent  swarm attack. 
\end{abstract}


\IEEEpeerreviewmaketitle

\section{Introduction}  

\IEEEPARstart{S}{warms} are characterized by large numbers of agents which act individually, yet produce collective, herd-like behaviors. Implementing cooperating swarm strategies for a large-scale swarm is a technical challenge which can be considered to be from the ``insider's perspective.'' It assumes inside control over the swarm's operating algorithms. However, as large-scale `swarm' systems of autonomous systems become achievable---such as those proposed by autonomous driving, UAV package delivery, and military applications---interactions with swarms outside our direct control becomes another challenge. This generates its own ``outsider's perspective'' issues. 

In this paper, we look at the specific challenge of protecting an asset against an adversarial swarm. Autonomous defensive agents are tasked with protected a High Value Unit (HVU) from an incoming swarm attack. The defenders do not fully know the cooperating  strategy employed by the adversarial swarm. Nevertheless, the task of the defenders is to maximize the probability of survival of the HVU against an attack by such a swarm. This challenge raises many issues---for instance, how to search for the swarm \cite{ifac2014}, how to observe and infer swarm operating algorithms \cite{JGCD2019}, and how to best defend against the swarm given algorithm unknowns and only limited, indirect control through external means. In this paper we restrict ourselves to the last issue. However, these problems share multiple technical challenges. The preliminary approach we apply in this paper demonstrates some basic methods which we hope will provoke development of more sophisticated tools. 

For objectives achieved via external control of the swarm, several features of swarm behavior must be characterized: capturing the dynamic nature of the swarm, tracking the collective risk profile created by a swarm, and engaging with a swarm via dynamic inputs such as autonomous defenders. The many modeling layers create a challenge for generating an effective response to the swarm, as model uncertainty and model error are almost certain. In this paper, we look at several dynamic systems where the network structure is determined by parameters. These parameters set neighborhood relations and interaction rules. Additional parameters establish defender input and swarm risk. 

We consider the generation of optimal defense strategies given uncertainty in parameter values. We demonstrate that small deviances in parameter values can have catastrophic effects on defense trajectories optimized without taking error into account. We then demonstrate the contrasting robustness of applying an uncertain parameter optimal control framework instead of optimizing with nominal values. The robustness against these parameter values suggests that refined parameter knowledge may not be necessary given appropriate computational tools. These computational tools---and the modeling of the high-dimensional swarm itself---are expensive. To assist with this issue, we provide dual conditions for this problem in the form of a Pontryagin minimum principle and prove the consistency of these conditions for the numerical algorithm. These dual conditions can thus be computed from the numerical solution of the computational method and provide a tool for solution verification and parameter sensitivity analysis.

The structure of this paper is as follows. Section \ref{modeling section}  provides examples of dynamic swarming models and extensions for defensive interactions.  
Section \ref{ProblemFormulation} discusses optimization challenges and describes  a general  uncertain parameter optimal control framework this problem could be addressed with. 
Section \ref{adjointconvergencesection} provides a proof of the consistency of the dual problem for this control framework, which expands on the results initially presented in the conference paper \cite{consistency_paper2}. Section \ref{numerical examples section} gives an example numerical implementation that demonstrates optimal defense against a large-scale swarm of $100$ agents. The final section discusses results and future work.

\section{Modeling Adverserial Swarms}\label{modeling section}  

\subsection{Cooperative Swarm Models} \label{swarm models}

The literature on the design of swarm strategies which produce coherent, stable collective behavior has become vast. A quick review of the literature points to two main trends/categories in swarm behavior design. The first one relies on dynamic modeling of the agents and potential functions to control their behavior (see \cite{TR2018,mehmood} and references therein).  The second trend uses rules to describe agents' motion and local rule-based algorithms to control them \cite{SwarmIntelligence}, \cite{SwarmIntelligence2}.

We present two examples of dynamic swarming strategy from the literature. These examples are illustrative of the forces considered in many swarming models:
\begin{itemize}
\item collision avoidance between swarm members
\item alignment forces between neighboring swarm members
\item stabilizing forces
\end{itemize}

These intra-swarm goals are aggregated to provide a swarm control law, which we will refer to as $F_S$, to each swarm agent. Both example models in this paper share the same double integrator form with respect to this control law. For $n$ swarm agents, dynamics are defined by
\begin{equation}
\begin{aligned}
& \ddot x_i = u_i. ~~ i = 1, \ldots, n,
\end{aligned}
\end{equation}
\begin{equation}
u_i = F_S(x_i, \dot{x}_i,  \forall j \neq i: x_j, \dot{x}_j | \theta).
\end{equation}

\subsubsection{Example Model 1: Virtual Body Artificial Potential, \cite{leonard2001virtual,ogren2004cooperative} }

In this model, swarm agents track to a virtual body (or bodies) guiding their course while also reacting to intra-swarm forces of collision avoidance and group cohesion. The input $u_i$ is the sum of intra-swarm forces, virtual body tracking, and a velocity dampening term. In addition, in this adversarial scenario, swarm agents are influenced to avoid intruding defense agents. The intra-swarm force between two swarm agents has magnitude $f_I$ and is a gradient  of an artificial potential $V_I$. Let 
\begin{equation}
x_{ij} = x_i - x_j.
\end{equation}
The artificial potential $V_I$ depends on the distance $||x_{ij}||$ between swarm agents $i$ and $j$. The artificial potential $V_I$ is defined as: 
\begin{equation}
V_I = \left\{ 
\begin{gathered}
  \alpha \left( \ln \left( ||x_{ij}|| \right)+\frac{d_0}{||x_{ij}||} \right), {~~~ 0 < ||x_{ij}|| < d_1} \hfill \\
  \alpha \left( \ln(d_1)+\frac{d_0}{d_1} \right), {~~~~~~~~~~~~~||x_{ij}|| \geq d_1} \hfill \\ 
\end{gathered}  \right. \label{eq:V_I}\\
\end{equation}
where $\alpha$ is a scalar control gain, $d_0$ and $d_1$ are scalar constants for distance ranges.
Then the magnitude of interaction force is given by
\begin{equation}
f_I = \left\{ 
\begin{gathered}
  \nabla_{||x_{ij}||} V_I, {~~~ 0 < ||x_{ij}|| < d_1} \hfill \\
  0, {~~~~~~~~~~~~~||x_{ij}|| \geq d_1} \hfill \\ 
\end{gathered}  \right. \label{eq:f_I}\\
\end{equation}
The swarm body is guided by `virtual leaders', non-corporeal reference trajectories which lead the swarm. We assign a potential $V_h$ on a given swarm agent $i$ associated with the $k$-th virtual leader, defined with the distance $||h_{ik}||$ between the swarm agent $i$ and leader $k$. Mirroring the parameters $\alpha$, $d_0$, and $d_1$ defining $V_I$, we assign $V_h$ the parameters $\alpha_h$, $h_0$, and $h_1$. An additional dissipative force $ f_{v_i}$ is included for stability. The control law $u_i$ for the vehicle $i$ associated with $m$ defenders is given by
\begin{equation}
\begin{aligned}
u_i & = -\sum_{j \neq i}^{n} \nabla _{x_i} V_I(x_{ij}) - \sum_{k = 1}^{m} \nabla _{x_i} V_h(h_{ik})  + f_{v_i} \\ 
&= -\sum_{j \neq i}^{n} \frac{f_I(x_{ij})}{||x_{ij}||} x_{ij}  - \sum_{k = 1}^{m} \frac{f_h(h_{ik})}{||h_{ik}||} h_{ik}  
 + f_{v_i}.
\end{aligned}
\label{eq:u_i}
\end{equation}

\subsubsection{Example Model 2:  Reynolds Boid Model, \cite{reynolds1987flocks}, \cite{mehmood} }
 For radius $r$, $j = 1, \dots, N$, define the neighbors of agent $i$ at position $x_i\in \mathbb{R}^n$by the set
  \begin{equation}
  \mathcal{N}_i = \{j | j\neq i \wedge  \|x_i-x_j \| < r\}
  \end{equation}
  Swarm control is designated by three forces.
  
{\it Alignment of velocity vectors:}
\begin{equation}
  f_{al} = -w_{al} \left( \dot{x}_i -  \frac{1}{| \mathcal{N}_i |} \sum_{j \in   \mathcal{N}_i } \dot{x}_j   \right)
 \end{equation} 
  
{\it Cohesion of swarm:}
 \begin{equation}
   f_{coh} = -w_{coh}  \left(  x_i - \frac{1}{| \mathcal{N}_i |} \sum_{j \in   \mathcal{N}_i } x_j   \right)
\end{equation}

{\it Separation between agents:}
 \begin{equation}\label{eqn: reynolds separation}
   f_{sep} = -w_{sep}  \frac{1}{| \mathcal{N}_i |}  \left(  \sum_{j \in   \mathcal{N}_i } \frac{x_j - x_i}{\| x_i - x_j \|} \right) 
   \end{equation}
 for positive constant parameters  $w_{al}$, $w_{coh}$, $w_{sep}$. 
  
 \begin{equation}
 u_i =   f_{al} + f_{coh} + f_{sep} 
 \end{equation}


 \subsection{Adversarial Swarm Models}\label{defense model}

In order to enable adversarial behavior and defense, the inner swarm cooperative forces $F_S$ need to be supplemented by additional forces of exogenoeous input into the collective. As written, the above cooperative swarming models neither respond to outside agents nor `attack' (swarm towards) a specific target. 

The review \cite{TR2018} discusses several approaches to adversarial control. Examples include containment strategies modeled after dolphins \cite{Egerstedt}, sheep-dogs \cite{sheepdogs1, sheepdogs2}, and birds of prey \cite{herdingbirds}. In \cite{Schwartz} the authors study interaction between two swarms, one of which can be considered adversarial. In these examples of adversarial swarm control, the mechanism of interaction and defense is provided through the swarm's own pursuit and evasion responses. This indirectly uses the swarm's own response strategy against it, an approach which can be termed `herding.'


We summarize a control scheme wuth HVU target tracking and herding driven by the reactive forces of collision avoidance with the defenders as the following, for HVU states $y_0$ and defender states $y_k$, $k=1, \dots, K$:
\begin{align*}
u_i =&F_S(x_i, \dot{x}_i,  \forall j \neq i: x_j, \dot{x}_j | \theta)   &{\color{red} \leftarrow \text{\it intra-swarm}} \\
+ &F_{HVU}(x_i, \dot{x}_i,  y_0, \dot{y}_0 | \theta)  &{\color{red} \leftarrow \text{\it target tracking}} \\
+ &F_D(x_i, \dot{x}_i,  \forall k \neq i: y_k, \dot{y}_k | \theta)  &{\color{red} \leftarrow \text{\it herding}} \numberthis \label{full swarm dynamics}
\end{align*}


In addition to herding reactions, one can consider more direct additional
forces of disruption, to model neutralizing swarm agents an/or physically removing them from the swarm. One form this can take, for example, is removal of agents from the communications network, as considered in \cite{szwaykowska2016collective}. Another approach is taken in \cite{walton2018optimal}, which uses survival probabilities based on damage attrition. Defenders and the attacking swarm engage in mutual damage attrition while the swarm also damages the HVU when in proximity. Probable damage between agents is tracked as damage rates over time, where the rate of damage is based on features such as distance between agents and angle of attack. The damage rate at time $t$ provides the probability of a successful `hit' in time period $[t, t+\Delta t]$. Probability of agent survival can be modeled based on the aggregate number of hits it takes to incapacitate. Reference \cite{walton2018optimal} provides derivations for multiple possibilities, such as single shot destruction and N shot destruction. These probabilities take the form of ODE equations. Tracking survival  probabilities thus adds an additional state to the dynamics of each agents---a survival probability state.

\subsubsection{Example Attrition Model: Single-Shot Destruction, \cite{walton2018optimal}}
Let $P_0(t)$ be the probability the HVU has survived up to time $t$, $P_k(t)$, $k=1, \dots, K$, the probability defender $k$ has survived, and $Q_j(t)$, $j=1, \dots, N$ the probability swarm attacker $j$ has survived.  Let $d_y^{j,k}(x_j(t),y_k(t))$ be the damage the defender $y_k$ inflicts on swarm attacker $x_j$ and let $d_x^{k,j}(y_k(t), x_j)$ be the damage the swarm attacker $x_j$ inflicts on the defender $y_k$, with the HVU represented by $k=0$. 

Then the survival probabilities for attackers and defenders from single shot destruction are given by the coupled ODEs:
\[
\begin{cases}
\dot{Q}_j(t)= & \\
-Q_j(t) \sum_{k=1}^K P_k(t) d_y^{j,k}(x_j(t),y_k(t)), & Q_j(0) = 1 \\
\dot{P}_k(t)= & \\
-P_k(t)\sum_{j=1}^N Q_j(t) d_x^{k,j}(y_k(t), x_j(t)), & P_k(0) = 1
\end{cases}
\]
for $j = 1, \dots, N$, $k=0\dots, K$.
  
\section{Problem Formulation}  
\label{ProblemFormulation}
\begin{figure}[h]
	\centering
		\includegraphics[scale=.23]{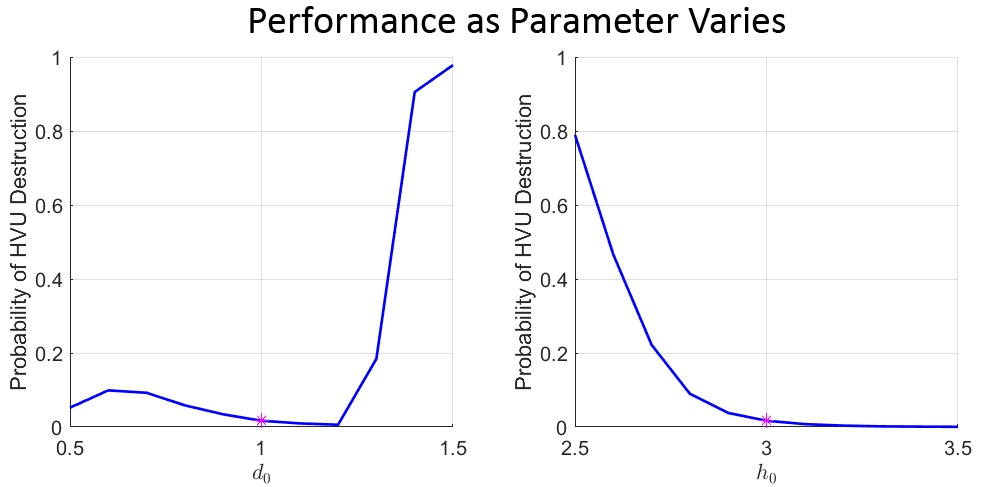}
		\caption{Example performance of solutions calculated using nominal values when parameter value is varied. Calculated using values in Section \ref{leonard_numerical_sec}. Magenta dot marks the nominal value used in optimization problem.}
	\label{nomex_fig}
\end{figure}

The above models depend on a large number of parameters. The dynamic swarming model coupled with attrition functions result in over a dozen key parameters, and many more would result from a non-homogeneous swarm. A concern would be that this adds too much model specificity, making optimal defense strategies lack robustness due to sensitivity to the specific set of model parameters. This concern turns out to be justified. When defense strategies are optimized for fixed, nominal parameter values, they display catastrophic failure for small perturbations of certain parameters as can be seen in Figure \ref{nomex_fig}.  In fact, the plots included in Figure \ref{nomex_fig} clearly demonstrate that the sensitivity of the cost with respect to the uncertain parameters is highly nonlinear. Thus, generating robust defense strategies requires a more sophisticated formalism introduced next in Section \ref{Uncertain Parameter Optimal Control}.

\subsection{Uncertain Parameter Optimal Control} \label{Uncertain Parameter Optimal Control}

The class of problems addressed by the computational algorithm is defined as follows:
\vspace{-5pt}
\\
\\ 
\noindent \textbf{Problem $\mathbf{P}$.}\ \ Determine the function pair $(x,u)$ with $x \in W_{1,\infty}([0,T]\times \Theta;\mathbb{R}^{n_x})$, $u \in L_\infty([0,T];\mathbb{R}^{n_u})$ that minimizes the cost
\small
\begin{align}\label{eqn:bobjective}
J[x,u] = \! \! \int_\Theta \! \Big[F\left(x(T,\theta),\theta\right) \!+\! \! \int_0^T\! r(x(t,\theta),u(t),t,\theta)dt\Big] d\theta  
\end{align}
\normalsize
subject to the dynamics 
\small
\begin{align}
\frac{\partial x}{\partial t}(t,\theta) & =  f(x(t,\theta),u(t),\theta), \label{eqn:dyn}
\end{align}
\normalsize
initial condition $x(0,\theta)  = x_0(\theta)$,
and the control constraint $g(u(t)) \leq 0$ for all $t \in [0,T]$. The set $L_\infty([0,T]; \mathbb{R}^{n_u})$ is the set of all essentially bounded functions, $W_{1,\infty}([0,T]\times \Theta;\mathbb{R}^{n_x})$ the Sobolev space of all essentially bounded functions with essentially bounded distributional derivatives,   and $F:\mathbb{R}^{n_x} \times \mathbb{R}^{n_\theta} \mapsto \mathbb{R}$, $r:\mathbb{R}^{n_x}\times\mathbb{R}^{n_u}\times\mathbb{R} \times \mathbb{R}^{n_\theta} \mapsto \mathbb{R}$, $g:\mathbb{R}^{n_u} \mapsto \mathbb{R}^{n_g}$. Additional conditions imposed on the state and control space and component functions are specified in Appendix Section \ref{regularity assumptions}.

In Problem $\mathbf{P}$, the set $\Theta$ is the domain of a parameter $\theta \in \mathbb{R}^{n_\theta}$.  The format of the cost functional is that of the integral over $\Theta$ of a Mayer-Bolza type cost with parameter $\theta$. This parameter can represent a range of values for a feature of the system, such as in ensemble control \cite{ensemble_inhomogeneous}, or a stochastic parameter with a known probability density function.


For computation of numerical solutions, we introduce an approximation of Problem $\mathbf{P}$, referred to as Problem $\mathbf{P^M}$. Problem $\mathbf{P^M}$ is created by approximating the parameter space, $\Theta$, with a numerical integration scheme. This numerical integration scheme is defined in terms of a finite set of $M$ nodes $\{\theta_i^M\}_{i=1}^{M}$ and an associated set of $M$ weights $\{\alpha_i^M\}_{i = 1}^M \subset \mathbb{R}$ such that
\begin{equation}\label{equ:quadrature}
\int_{\Theta} h(\theta) d\theta  =  \lim_{M \to \infty} \sum_{i = 1}^M h(\theta^M_i) \alpha_i^M.
\end{equation}
given certain function smoothness assumptions. See Appendix Assumption \ref{numericalschemeassumption} for formal assumptions. Throughout the paper, $M$ is used to denote the number of nodes used in this approximation of parameter space. 
 
For a given set of nodes $\{\theta_i^M\}_{i=1}^{M}$, and control $u(t)$, let $\bar{x}_i^M(t)$, $i=1,\dots,M$, be defined as the solution to the ODE created by the state dynamics of Problem $\mathbf{P}$ evaluated at $\theta_i^M$:
\small
\begin{align}\label{eqn:barxdyn}
\begin{cases}
\frac{d\bar{x}_i^M}{dt}(t) = f( \bar{x}_i^M(t),u(t),\theta_i^M) \\
 \bar{x}_i^M(0) = x_0(\theta_i^M),
\end{cases}
\quad i=1,\dots,M.
\end{align}
\normalsize
Let $\bar{X}^M(t) = [ \bar{x}_1^M(t),\ldots, \bar{x}_M^M(t)]$. The system of ODEs defining $\bar{X}^M$ has dimension $n_x\times M$, where $n_x$ is the dimension of the original state space and $M$ is the number of nodes. The numerical integration scheme for parameter space creates an approximate objective functional, defined by:
\small
\[
\bar{J}^M\hspace{-1pt} [ \bar{X}^M\hspace{-1pt} ,u] \hspace{-2pt}  = \qquad \qquad  \qquad \qquad  \qquad \qquad  \qquad
\]
\begin{equation}\label{eqn:jmobjective}
  \quad \sum_{i = 1}^M \Big[  F\left( \bar{x}_i^M (T),\theta_i^M \right)  +\int^T_0  r( \bar{x}_i^M (t),u(t),t,\theta_i^M )dt\hspace{-1pt}  \Big] \alpha_i^M.
 \end{equation}
\normalsize



In  \cite{walton_IJC}, the consistency of $\mathbf{P^M}$ is proved. This is the property that if optimal solutions to Problem $\mathbf{P^M}$ converge as the number of nodes $M\rightarrow \infty$, they converge to feasible, optimal solutions of Problem $\mathbf{P}$. See \cite{walton_IJC} for detailed proof and assumptions.

\subsection{Computational Efficiency}

The computation time of the numerical solution to the discretized problem defined in equations (\ref{eqn:barxdyn},\ref{eqn:jmobjective}) will depend on the value of $M$. Ideally, it should be sufficiently small as to allow for a fast solution. On the other hand, a value of $M$ that is too small will result in a solution that is not particularly useful, i.e too far from the optimal. Naturally, the question arises: how far is a particular solution from the optimal? One tool for assessing this lies in computing the Hamiltonian and is addressed in the next section.

\section{Consistency of Dual Variables} \label{adjointconvergencesection}


The dual variables provide a method to determine the solution of an optimal control problem or a tool to validate a numerically computed solution.  For numerical schemes based on direct discretization of the control problem, analyzing the properties of the dual variables annd their resultant Hamiltonian may also lead to insight into the validity of approximation scheme \cite{Hager.00, Qi_covector1}. This could be especially helpful in high-dimensional problems such as swarming, where parsimonious discretization is crucial to computational tractability. 

\begin{figure}[h]
	\centering
		\includegraphics[scale=.42]{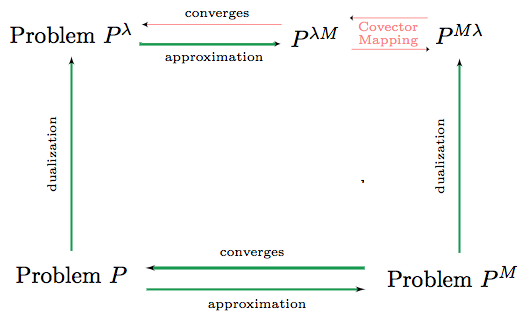}
		\caption{Diagram of primal and dual relations for parameter uncertainty control. Red lines designate the contribution of this paper.}
	\label{covector_partial}
\end{figure}

Previous work shows the consistency of the primal variables in approximate Problem $\mathbf{P^M}$ to the original parameter uncertainty framework of Problem $\mathbf{P}$. Here we build on that and prove the consistency of the dual problem of Problem $\mathbf{P}$ as well. This theoretical contribution is diagrammed in Figure \ref{covector_partial}. The consistency of the dual problem in parameter space enables approximate computation of the Hamiltonian from numerical solutions.  

In \cite{gabasov} necessary conditions for Problem $\mathbf{P}$ (subject to the assumptions of Section \ref{regularity assumptions}) were established. These conditions are as follows: 
\\
\\ 
\noindent \textbf{Problem $\mathbf{P^{\lambda}}$.}\quad\textbf{(\cite{gabasov}, pp. 80-82)}.  If $(x^*,u^*)$ is an optimal solution to Problem $\mathbf{P}$, then there exists an absolutely continuous  costate vector $\lambda^*(t,\theta)$ such that for $\theta \in \Theta$:
\normalsize 
\[
\frac{\partial \lambda^*}{\partial t}(t,\theta)  = -\frac{\partial H(x^*,\lambda^*,u^*,t,\theta)}{\partial x}
,
\]
\begin{equation}\label{adjoint}
 \lambda^*(T,\theta) = \frac{\partial F(x^*(T,\theta),\theta) }{\partial x} 
\end{equation}
where $H$ is defined as:
\[ 
H(x,\lambda,u,t,\theta) = \qquad  \qquad  \qquad  \qquad  \qquad  \qquad 
\]
\begin{equation}\label{first_Hamiltonian}
\lambda f(x(t,\theta),u(t),\theta) + r(x(t,\theta), u(t),t,\theta). 
\end{equation}
Furthermore, the optimal control $u^*$ satisfies
\begin{equation*}
u^*(t) = \underset{u \in U}{\arg\min} \mathbf{H}(x^*,\lambda^*,u,t),
\end{equation*} 
where $\mathbf{H}$ is given by
\begin{equation}\label{second_Hamiltonian_def}
\mathbf{H}(x,\lambda,u,t)= \int_{\Theta} H(x,\lambda, u,t,\theta) d\theta.
\end{equation}


Because Problem $\mathbf{P^M}$ is a standard nonlinear optimal control problem, it admits a dual problem as well, Problem $\mathbf{P^{M\lambda}}$, provided by the Pontryagin Minimum Principle (a survey of Minimum Principle conditions is given by\cite{Hartletall.95}). Applied to  $\mathbf{P^M}$ this generates:
\\
\\
{\bf Problem ${\mathbf{P^{M \lambda}}}$:} For feasible solution $(\bar{X}^M,u)$ to Problem ${\mathbf{P^M}}$, find $\bar{\Lambda}(t)=[\bar{\lambda}_1^M(t),  \dots \bar{\lambda}_M^M(t)]$, $\bar{\lambda}_i^M:[0,T] \rightarrow \mathbb{R}^{N_x}$, that satisfies the following conditions:
\[
\frac{d\bar{\lambda}_i^M}{dt}(t)  = -\frac{\partial \bar{H}^M(\bar{x}_i^M,\bar{\lambda}_i^M,u,t)}{\partial x_i^M}
,
\]
\begin{equation}\label{PMlambda adjoint}
 \bar{\lambda}_i^M(T) = \alpha_i^M \frac{\partial F(\bar{x}_i^M,\theta_i^M) }{\partial \bar{x}_i^M}, 
\end{equation}
where $\bar{H}^M$ is defined as:
\[
\bar{H}^M(\bar{X}^M,\bar{\Lambda}_M,u,t) =
\qquad \qquad \qquad \qquad \qquad \qquad \qquad 
\]
\[
\sum_{i=1}^M \left[ \bar{\lambda}_i^M f(\bar{x}_i^M(t),u(t),\theta_i^M) + \right. \qquad \qquad 
\]
\begin{equation}\label{PMlambda Hamiltonian}
\qquad \qquad \left. \alpha_i^M r(\bar{x}_i^M(t), u(t),t,\theta_i^M) \right].
\end{equation}

An alternate direction from which to approach to solving Problem $\mathbf{P}$ overall is to approximate the necessary conditions of Problem $\mathbf{P}$ , i.e. Problem $\mathbf{P^\lambda}$, directly rather than to approximate Problem $\mathbf{P}$. This creates the system of equations:
\[
\frac{d\lambda}{dt}(t,\theta_i^M)  = -\frac{\partial H(x,u,t,\theta_i^M)}{\partial x} 
\]
\begin{equation}\label{collocated lambdas}
\lambda(T,\theta_i^M) = \frac{\partial F(x(T,\theta_i^M),\theta_i^M) }{\partial x} 
\end{equation}
for $i=1,\dots,M$, where $H$ is defined as:
\[
H(x,\lambda,u,t,\theta) =  \qquad \qquad  \qquad \qquad  \qquad \qquad 
\]
\[
 \qquad \qquad  \lambda f(x(t,\theta),u(t),\theta) + r(x(t,\theta), u(t),t,\theta). 
\]
This system of equations can be re-written in terms of the quadrature approximation of the stationary Hamiltonian defined in equation (\ref{second_Hamiltonian_def}). Define
\[
\tilde{H}^M(x,\lambda,u,t): = \qquad \qquad  \qquad \qquad  \qquad \qquad 
\]
\[
\sum_{i=1}^M \alpha_i^M H(x(t,\theta_i^M),\lambda(t,\theta_i^M), u(t),t,\theta_i^M)
.
\]
Let
\[\tilde{\Lambda}(t)
=[\tilde{\lambda}_1^M(t),  \dots \tilde{\lambda}_M^M(t)]
=[\lambda(t,\theta_1^M),\dots, \lambda(t,\theta_M^M)]
\]
and let
\[
\tilde{X}_M = [\tilde{x}_1^M(t),\dots,\tilde{x}_M^M(t)]
\] 
denote the semi-discretized states from equation (\ref{eqn:barxdyn}). Equation (\ref{collocated lambdas}) can then be written as: 
\[
\frac{d \tilde{\lambda}_i^M}{dt}(t) = -\frac{1}{\alpha_i^M}\cdot\frac{\partial \tilde{H}^M(\tilde{X}_M,\tilde{\Lambda},u,t)}{\partial \tilde{x}_i^M} 
\]
\begin{equation} 
\tilde{\lambda}_i^M(T) = \frac{\partial F(\tilde{x}_i^M(T), \theta_i^M) }{\partial \tilde{x}_i^M}  
\end{equation}
for $i=1,\dots,M$. Thus we reach the following discretized dual problem:
\\
\\
{\bf Problem ${\mathbf{P^{ \lambda M}}}$:} For feasible contols $u$ and solutions $\tilde{X}_M$ to equation (\ref{eqn:barxdyn}), find $\tilde{\Lambda}(t)=[\tilde{\lambda}_1^M(t),  \dots \tilde{\lambda}_M^M(t)]$, $\tilde{\lambda}_i^M:[0,T] \rightarrow \mathbb{R}^{n_x}$, that satisfies the following conditions:
\[
\frac{d\tilde{\lambda}_i^M}{dt}(t) = -\frac{1}{\alpha_i^M}\cdot\frac{\partial \tilde{H}^M(\tilde{X}_M,\tilde{\Lambda},u,t)}{\partial \tilde{x}_i^M}
, 
\]
\begin{equation}\label{PlambdaM adjoint}
\tilde{\lambda}_i^M(T) = \frac{\partial F(\tilde{x}_i^M,\theta_i^M) }{\partial \tilde{x}_i^M}, 
\end{equation}
where $\tilde{H}^M $  is defined as:
\[
\tilde{H}^M(\tilde{X}_M,\tilde{\Lambda}_M,u,t) =
\qquad \qquad \qquad \qquad \qquad \qquad \qquad 
\]
\[
\sum_{i=1}^M \left[ \alpha_i^M \tilde{\lambda}_i^M f(\tilde{x}_i^M(t),u(t),\theta_i^M) + \right.
\qquad \qquad \qquad  
\]
\begin{equation}\label{PlambdaM Hamiltonian}
\qquad \qquad \qquad \qquad \qquad 
\left.  \alpha_i^M r(\tilde{x}_i^M(t), u(t),t,\theta_i^M) \right].
\end{equation}

In the case of this particular problem, unlike standard control,  the collocation of the relevant dynamics involves no approximation of differentiation (since the discretization is in the parameter domain rather than the time domain) and thus the mapping of covectors between Problem ${\mathbf{P^{M \lambda }}}$ and Problem ${\mathbf{P^{ \lambda M}}}$ is straightforward. 
\begin{lemma}\label{covector lemma} 
The mapping:
\[
(\bar{x}_i^M,\bar{u}) \mapsto (\tilde{x}_i^M,\tilde{u} ), \quad
\frac{\bar{\lambda}_i^M}{\alpha_i^M} \mapsto \tilde{\lambda}_i^M, 
\]
for $i=1, \dots, M$ is a bijective mapping from solutions of Problem ${\mathbf{P^{M \lambda }}}$ to Problem ${\mathbf{P^{ \lambda M}}}$.
\end{lemma}

\begin{theorem}\label{hamiltonian convergence theorems}
Let $\{\tilde{X}_M, \tilde{\Lambda}_M, u_M\}_{M\in V}$ be a sequence of solutions for Problem $\mathbf{P^{\lambda M}}$ with an accumulation point $\{\tilde{X}^\infty, \tilde{\Lambda}^\infty, u^\infty\}$. Let $(x^\infty,\lambda^\infty,u^\infty)$ be the solutions to Problem $\mathbf{P^{\lambda }}$ for the control $u^\infty$. Then
\[
\lim_{M\in V} \tilde{H}^M(\tilde{X}_M, \tilde{\Lambda}_M, u_M,t)
=
\mathbf{H}(x^\infty,\lambda^\infty,u^\infty,t)
\]
where $\tilde{H}^M$ is the Hamiltonian of Problem $\mathbf{P^{\lambda M}}$ as defined by equation (\ref{PlambdaM Hamiltonian}) and $\mathbf{H}$ is the Hamiltonian of Problem $\mathbf{P}$ as defined by equation (\ref{second_Hamiltonian_def}). The proof of this theorem can be found in the Appendix.
\end{theorem}

The convergence of the Hamiltonians of the approximate, standard control problems to the Hamiltonian of the general problem, $\mathbf{H}(x^\infty,\lambda^\infty,u^\infty,t)$, means that many of the useful features of the Hamiltonians of standard optimal control problems are preserved. For instance, it is straightforward to show that 
the satisfaction of Pontryagin's Minimum Principle by the approximate Hamiltonians
implies minimization of $\mathbf{H}(x^\infty,\lambda^\infty,u^\infty,t)$ as well. That is, that
\[
\mathbf{H}(x^\infty, \lambda^\infty,u^\infty,t) \leq  \mathbf{H}(x^\infty, \lambda^\infty,u,t)
\]
for all feasible $u$. Furthermore, when applicable, the stationarity properties of the standard control Hamiltonian--such as a constant-valued Hamiltonian in time-invariant problems, or stationarity with respect to $u(t)$ in problems with open control regions--are also preserved.

\section{Numerical Example} \label{numerical examples section}



In a slight refashioning of the notation in the Section \ref{defense model}, equation (\ref{full swarm dynamics}), let the parameter vector $\theta$ be defined by all the {\it unknown} parameters defining the interaction functions. Assuming prior distribution $\phi(\theta)$ over these unknowns and parameter bounds $\Theta$, we construct the following optimal control problem for robustness against the unknown parameters.
\\
\\ 
\noindent{\bf Problem SD (Swarm Defense):} For $K$ defenders and $N$ attackers, determine the defender controls $u_k(t)$ that minimize:
\begin{equation}
 J =\int_\theta \left[1-P_0(t_f,\theta) \right] \phi(\theta)d\theta \label{SwarmDefense}
\end{equation}
subject to:
\small

\[
\begin{cases}
\dot{y}_k(t) = f(y_k(t), u_k(t)), & \hspace{-25pt} y_k(0) =  y_{k0}\\
\ddot{x}_j(t,\theta) = & \\
F_S(t,\theta) + F_{HVU}(t,\theta)+F_D(t,\theta), & \hspace{-25pt} x_j(0,\theta) = x_{j0}(\theta)  \\
\dot{Q}_j(t,\theta)= & \\
-Q_j(t,\theta) \sum_{k=1}^K P_k(t,\theta) d_y^{j,k}(x_j(t,\theta),y_k(t)), & \hspace{-5pt} Q_j(0,\theta) = 1 \\
\dot{P}_k(t)= & \\
-P_k(t,\theta)\sum_{j=1}^N Q_j(t,\theta) d_x^{k,j}(y_k(t), x_j(t,\theta)), & \hspace{-5pt}P_k(0,\theta) = 1
\end{cases}
\]
\normalsize
for swarm attackers $j = 1, \dots, N$ and controlled defenders $k=1\dots, K$.

We implement Problem {\bf SD} for both swarm models in Section \ref{swarm models}, for a swarm of $N=100$ attackers and $K=10$ defenders.

\subsection{Example Model 1: Virtual Body Artificial Potential} \label{leonard_numerical_sec}

The cooperative swarm forces $F_S$ are defined with the Virtual Body Artificial Potential of Section \ref{swarm models} with parameters $\alpha$, $d_0$ and $d_1$. In lieu of a potential for the virtual leaders, we assign the HVU tracking function:
\begin{equation}\label{eqn:hvu tracking}
f_{HVU} = - \frac{K_1(x_i-y_0)}{\|x_i-y_0\|}
\end{equation}
where $y_0\in\mathbb{R}^3$ is the position of the HVU. The dissipative force $f_{v_i} = - K_2\dot x_i$ is employed to guarantee stability of the swarm system. $K_1$ and $K_2$ are positive constants. The swarm's collision avoidance response to the defenders is defined by equation (\ref{eq:V_I}) with parameters $\alpha_h$, $h_0$ and $h_1$. Since there is only a repulsive force between swarm members and defenders, not an attractive force, we set $h_1=h_0$. For attrition, we use the the damage function defined in equation (21) of \cite{walton2018optimal}. For the damage rate of defenders inflicted on attackers, we calibrate by the parameters $\lambda_D$, $\sigma_D$. For the damage rate of attackers inflicted on defenders, we calibrate by the parameters $\lambda_A$, $\sigma_A$. In both cases, the parameters $F$ and $a$ in \cite{walton2018optimal} are set to $F=0$, $a=1$. Table \ref{Fixed Parameter Values1} provides the parameter values that remain fixed in each simulation, and and Table \ref{tbl:uncert-param-values1} provides the parameters we consider as uncertain.

\begin{table}[h!]
\centering
\begin{tabular}{c||  r |c|}
 { Parameter} & { Value} & { Reference}\\
 \hline
  $t_f$ & $45$ & {\em  final time}\\
\hline
  $K_1$ & $5$ & {\em  tracking coefficient}\\
\hline
  $K$ & $10$ & {\em  number of defenders}\\
\hline
  $h_1$ & $h_0$ & {\em  interaction parameter}\\
 \hline
 $\lambda_D$ & $2$ & {\em defender weapon intensity}\\
\hline
  $\sigma_D$  &  $2$ & {\em  defender weapon range}\\
\hline
  $N$ & $100$ & {\em  number of attackers}\\
\hline
   $K_2$ &  $5$ & {\em  dissipative force}\\
\hline

\noalign{\smallskip}
\end{tabular}
\caption{Model 1 Fixed Parameter Values}
\label{Fixed Parameter Values1}
\end{table}
\begin{table}[h!]
\centering
\begin{tabular}{c|| r |c|c|}
 { Parameter} & { Nominal} & { Range} & {Reference}\\
\hline
 $\alpha$ &  $.5$ & [0.1, 0.9] & {\em control gain}\\
\hline
  $d_0$  & $1$ & [0.5, 1.5] & {\em  lower range limit}\\
\hline
  $d_1$ & $6$ & [4, 8] & {\em upper range limit}\\
\hline
  $\lambda_A$ & $.05$ & [.01,.09] & {\em  weapon intensity}\\
\hline
   $\sigma_A$ & $2$  & [1.5 2.5] & {\em  weapon range}\\
\hline
$\alpha_h$ &  $6$ & [5,7] & {\em herding intensity} \\
\hline
  $h_0$  & $3$ & [2,4] & {\em herding range} \\
\hline
\noalign{\smallskip}
\end{tabular}
\caption{Model 1 Varied Parameter Values}
\label{tbl:uncert-param-values1}
\end{table}

\begin{figure}[h]
	\centering
		\includegraphics[scale=.48]{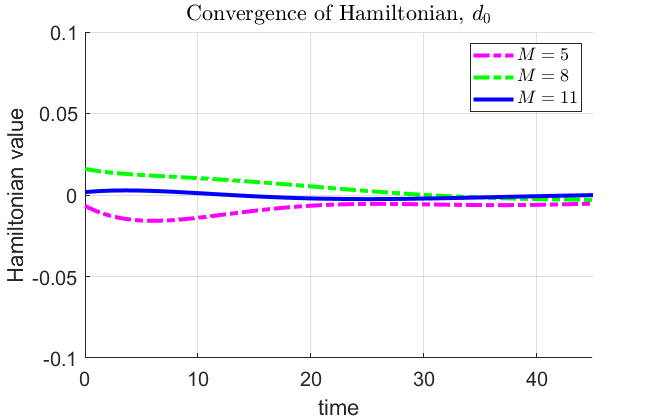}
		\caption{Convergence of Hamiltonion as number of parameter nodes $M$ increases}
	\label{leonard_H_fig}
\end{figure}

We first use the nominal parameter values provided in Tables \ref{Fixed Parameter Values1} and \ref{tbl:uncert-param-values1} to find a nominal solution defender trajectories that result in the minimum probability of HVU destruction. With the results of these simulations as a reference point, we consider as uncertain each of the parameters that define attacker swarm model and weapon capabilities. In this simulation, these parameters are considered individually. The number of discretization nodes for parameter space was chosen by examination of the Hamiltonian. To illustrate this method and the results obtained in Section \ref{adjointconvergencesection} we compute Hamiltonians for the Problem {\bf SD} and Model 1 with $\theta = d_0, d_0 \in [0.5,1.5]$ and $M = [5,8,11]$. As $M$ increases the sequence of Hamiltonians  should converge to the optimal Hamiltonian for the Problem  {\bf SD}. For Problem {\bf SD} that should result in a constant, zero-valued Hamiltonian. Figure \ref{leonard_H_fig} shows the respective Hamiltonians for $M = [5,8,11]$. The value $M=11$ was chosen for simulations, based on the approximately zero-valued Hamiltonian it generates.

\begin{figure}[ht] 
	\centering 
	\includegraphics[width = 0.49\columnwidth]{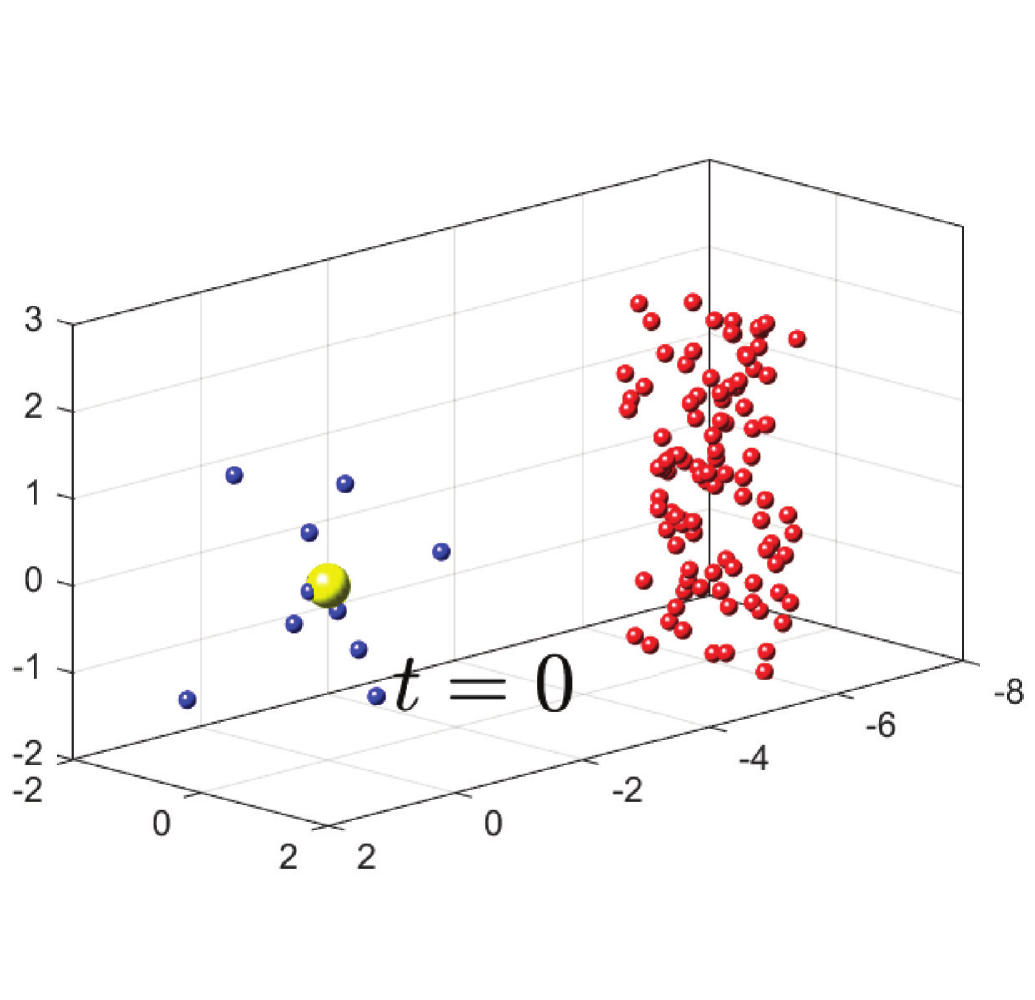}
	\includegraphics[width = 0.49\columnwidth]{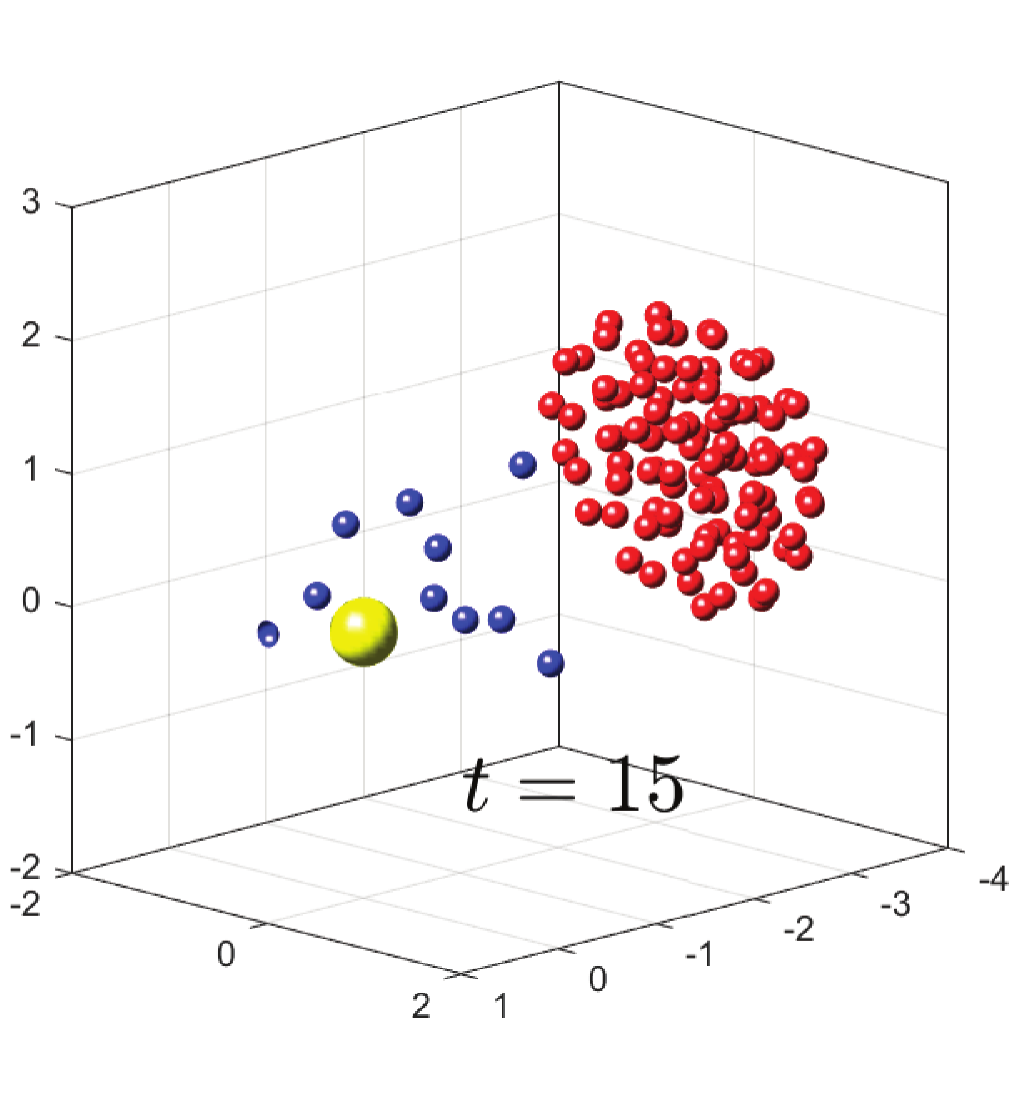}
	\includegraphics[width = 0.49\columnwidth]{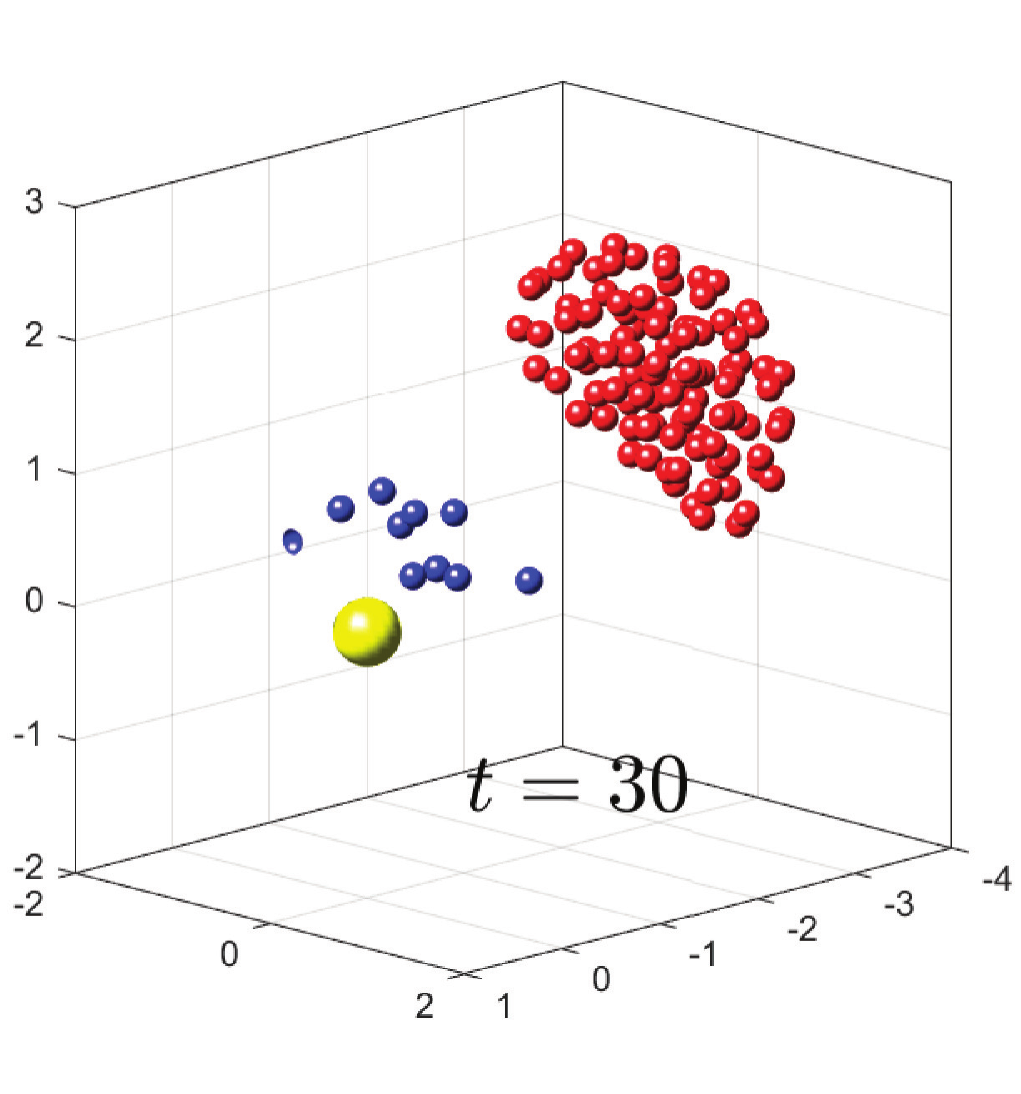}
	\includegraphics[width = 0.49\columnwidth]{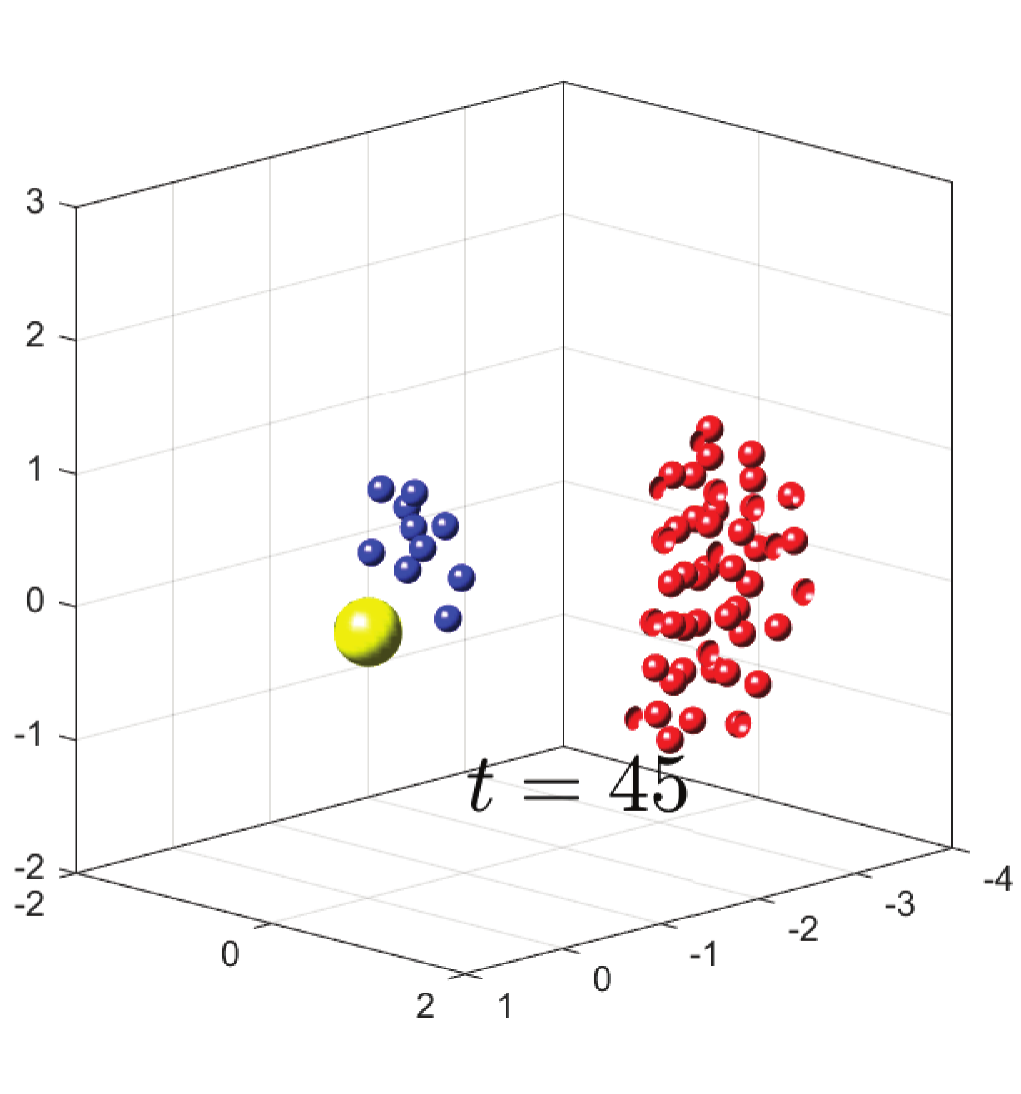}
	\includegraphics[width = \columnwidth]{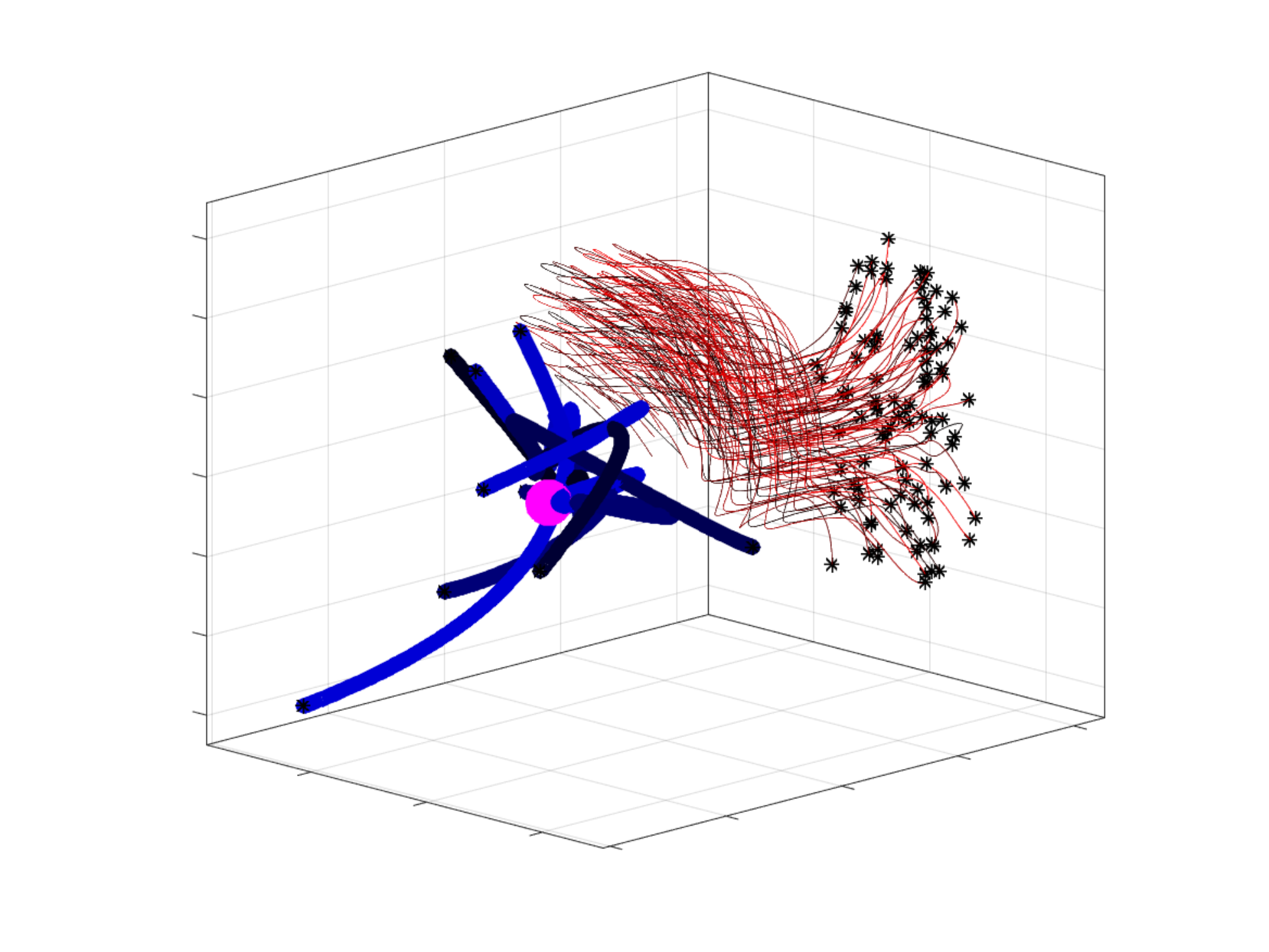}
	\caption{Shown are four snapshots during a simulations at $t=0$, 15, 30, and 45 (time units are arbitrary). Defenders are represented by blue spheres and attackers by red spheres. Below these snapshots, we show full trajectories for the entire simulation, which is the result of an optimization protocol using the parameters shown in Table~\ref{Fixed Parameter Values1} }.
	\label{fig:snapshots}
\end{figure}

\begin{figure*}[h]
	\centering
		\includegraphics[scale=.33]{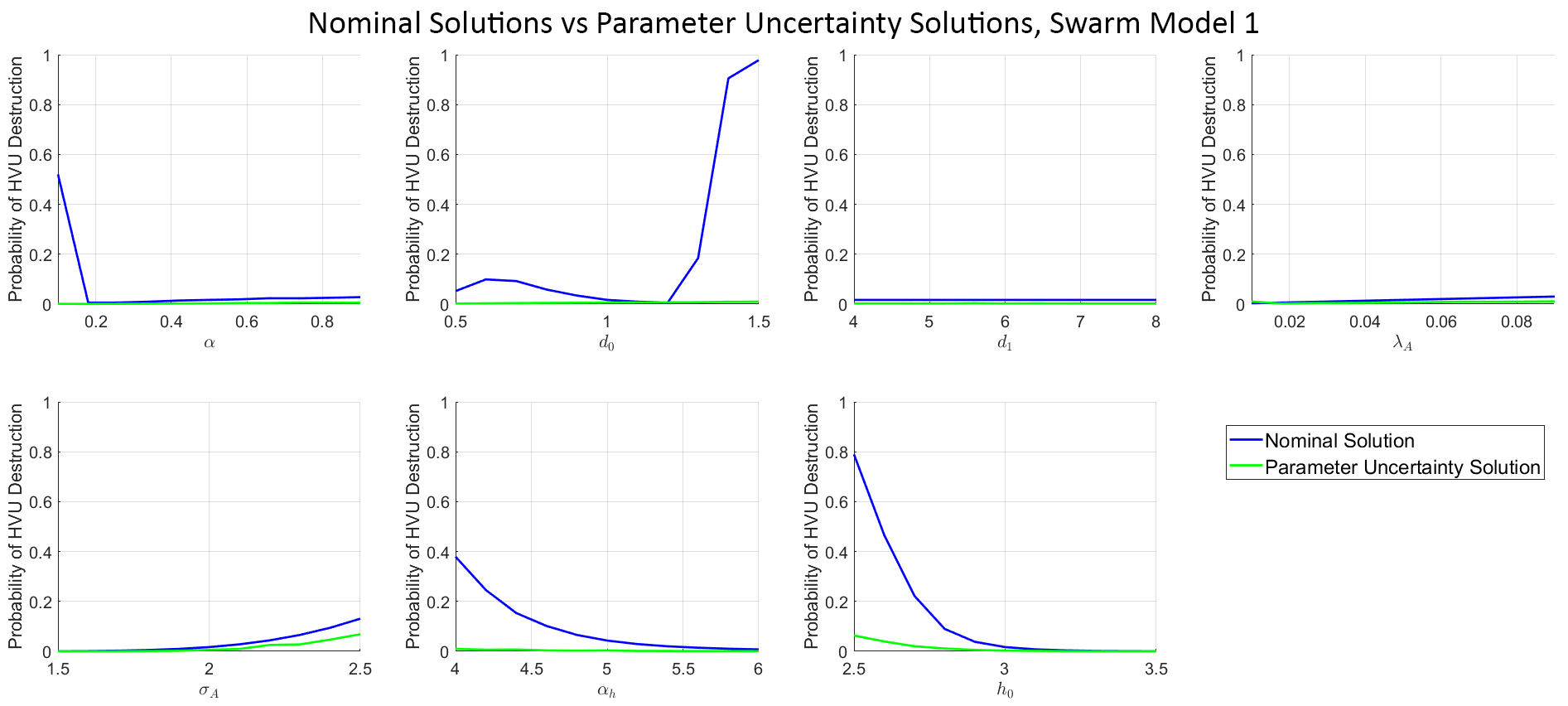}
		\caption{Performance of Solutions of Swarm Model 1 as parameter values are varied}
	\label{leonard_comparison_fig}
\end{figure*}

We compare the performance of the solution generated using uncertain parameter optimal control Problem {\bf SD} versus a solution obtained with the nominal values. Figure \ref{fig:snapshots} shows the nominal solution trajectories. The comparitive results of the nominal solutions vs the uncertain parameter control solutions are shown in Figure~\ref{leonard_comparison_fig}, where the performance of each is shown for different parameters values. 

As seen in Figure~\ref{leonard_comparison_fig} the trajectories generated by optimization using the nominal values perform poorly over a range of $\alpha$, $d_0$,   $\sigma_A$, $\alpha_k$ and $h_0$. In the case of $h_0$, for example, this is because the attackers are less repelled by the defenders when $h_0$ is decreased, and they are more able to destroy the HVU from a longer distance as $\sigma_A$ is increased. The parameter uncertainty solution, however, demonstrates that using the uncertain parameter optimal control framework a solution can be provided which is robust over a range of parameter values. We contrast these results with the case of uncertain parameters $d_1$ and $\lambda_A$, also shown in Figure~\ref{leonard_comparison_fig}.  It can be seen that robustness improvements are modest to non-existent for these parameters. This suggests an insensitivity of the problem $d_1$ and $\lambda_A$ parameters. This kind of analysis can be used to guide inference and observability priorities. 

\subsection{Example Model 2:  Reynolds Boid Model}

\begin{figure*}[h]
	\centering
		\includegraphics[scale=.33]{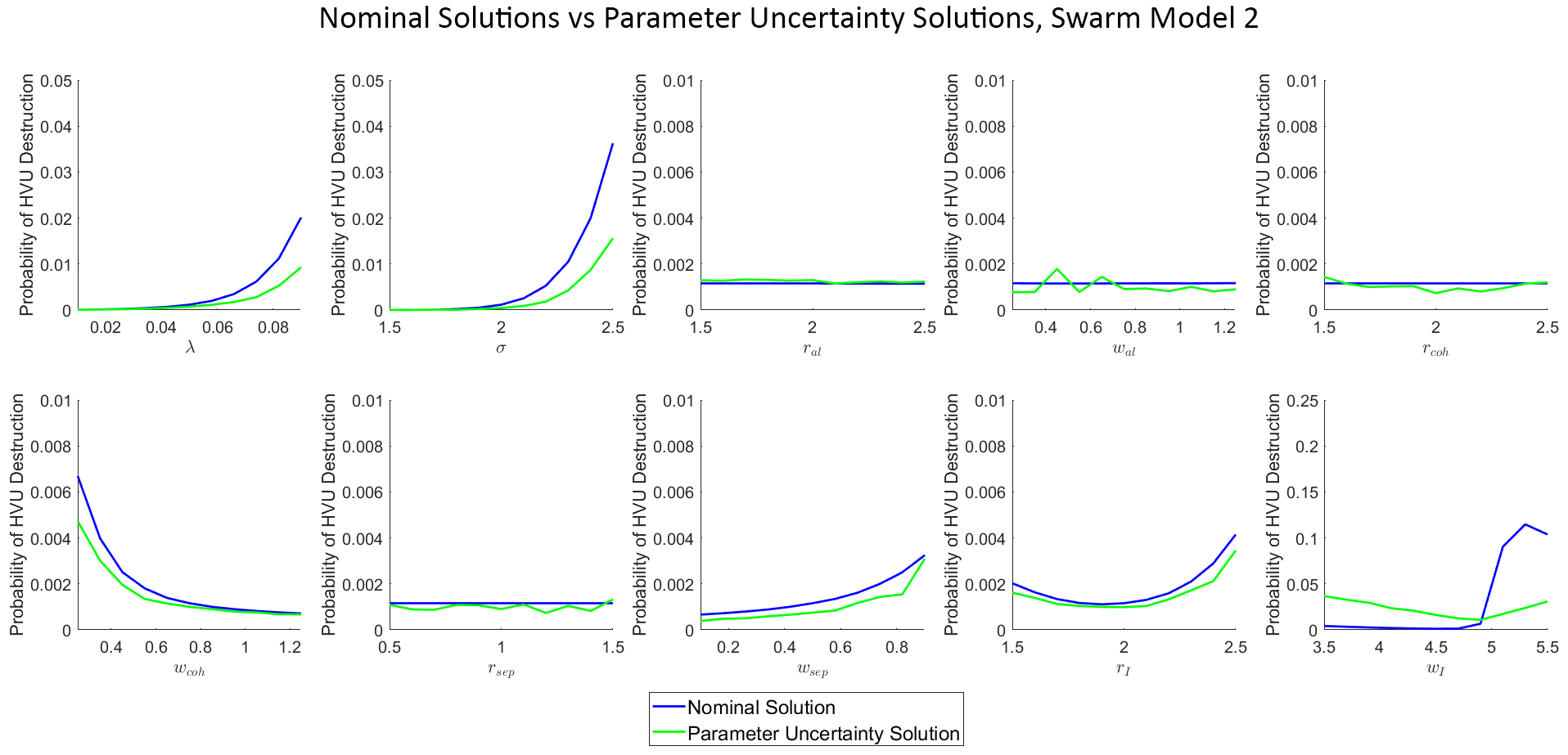}
		\caption{Performance of Solutions of Swarm Model 2 as parameter values are varied}
	\label{reynolds_comparison_fig}
\end{figure*}

To demonstrate flexibility of the proposed framework to include diverse swarm models we have applied the same analysis  as was done in Section  \ref{leonard_numerical_sec} to the Reynolds Boid Model introduced in Section \ref{swarm models}. We apply the same HVU tracking function as equation (\ref{eqn:hvu tracking}). The herding force $F_D$ of the defenders repelling attackers is applied as a separation force in the form of equation (\ref{eqn: reynolds separation}). The fixed parameter values are the same as those in Table \ref{Fixed Parameter Values1}; the uncertain parameters and ranges are given in Table \ref{tbl:uncert-param-values2}. The results are shown in Figure \ref{reynolds_comparison_fig}. Again, we see that the tools developed in this paper can be used to gain an insight into the robustness properties of the nominal versus uncertain parameter solutions. For example, we can see that the uncertain parameter solutions perform much better than the nominal ones for the cases where $\lambda$, $\sigma$ and $w_I$ are uncertain. 

\begin{table}[h!]
\centering
\begin{tabular}{c|| r |c|c|}
 { Parameter} & { Nominal} & { Range} & {reference}\\
\hline
  $\lambda_A$ & $.05$ & [.01,.09] & {\em  weapon intensity}\\
\hline
   $\sigma_A$ & $2$  & [1.5 2.5] & {\em  weapon range}\\
\hline

$r_{al}$		& $2$  & [1.5,2.5] & {\em alignment range}\\
$w_{al}$		& $.75$  & [.25,1.25] & {\em alignment intensity}\\
$r_{coh}$		& $2$  & [1.5,2.5] & {\em cohesion range}\\
$w_{coh}$	& $.75$  & [.25,1.25] & {\em cohesion intensity}\\
$r_{sep}$		& $1$  & [.5, 1.5] & {\em separation range}\\
$w_{sep}$	& $0.5$  & [.1, .9] & {\em separation intensity}\\
$r_{I}$		& $2$  & [1.5,2.5] & {\em herding range}\\
$w_{I}$		& $4.5$  & [3.5, 5.5] & {\em herding intensity}\\
\noalign{\smallskip}
\end{tabular}
\caption{Model 2 Varied Parameter Values}
\label{tbl:uncert-param-values2}
\end{table}

\section{Conclusions}
In this paper we have built on our previous work on developing an efficient  numerical framework for solving uncertain parameter optimal control problems.  Unlike uncertainties introduced into systems due to stochastic ``noise,'' parameter uncertainties do not average or cancel out in regards to their effects. Instead, each possible parameter value creates a specific profile of possibility and risk.  The uncertain optimal control framework which has been developed for these problems exploits this inherent structure by producing answers which have been optimized over all parameter profiles. This approach takes into account the possible performance ranges due to uncertainty, while also utilizing what information is known about the uncertain features--such as parameter domains and prior probability distributions over the parameters. Thus we are able to contain risk while providing plans which have been optimized for performance under all known conditions. The results reported in this paper include analysis of the consistency of the adjoint variables of the numerical solution. In addition, the  paper includes  a numerical analysis of a large scale adversarial swarm engagement  that clearly demonstrates the benefits of using the  proposed framework.

There are many directions of future work for the topics of this paper. The numerical simulations in this paper consider the parameters individually, as one-dimensional parameter spaces. However, Problem {\bf P} allows for multi-dimensional parameter spaces. A more dedicated implementation, taking advantage of the parallelizable form of equation (\ref{eqn:barxdyn}) for example, could certainly manage several simultaneous parameters. Exponential growth as parameter space dimension increases is an issue for both the quadrature format of equation (\ref{equ:quadrature}) and handling of the state space size for equation (\ref{eqn:barxdyn}). This can be somewhat mitigated by using sparse grid methods for high-dimensional integration to define the nodes in equation (\ref{equ:quadrature}). For large enough size, Monte Carlo sampling, rather than quadrature might be more appropriate for designating parameter nodes.

Another direction of future work is in greater application of the duality results of Section \ref{adjointconvergencesection}. The numerical results in this paper simply utilize the Hamiltonian consistency. The proof of Theorem \ref{hamiltonian convergence theorems}, however, additionally demonstrates the consistency of the adjoint variables for the problem. As the results demonstrate, parameter sensitivity for these swarm models is highly nonlinear. The numerical solutions of Section \ref{numerical examples section} are able to demonstrate this sensitivity by applying the solution to varied parameter values. However, this is actually a fairly expensive method for a large swarm, as it involves re-evaluation of the swarm ODE for each parameter value. More importantly, it would not be scalable to high-dimensional parameter spaces, as the exponential growth of that approach to sensitivity analysis would be unavoidable. The development of an analytical adjoint sensitivity method for this problem could be of great utility for paring down numerical simulations to only focus on the parameters most relevant to success.

\appendix[Proof of Theorem \ref{hamiltonian convergence theorems}]

\subsection{Assumptions and Definitions}\label{regularity assumptions}

We we impose the assumptions of \cite{walton_IJC} Section 2. The definition of accumulation point used in the following proof can be found in \cite{walton_IJC} Definition 3.2. The following assumption is placed on the choice of numerical integration scheme to be utilized in approximating Problem $\mathbf{P}$:
\begin{assumption}\label{numericalschemeassumption}
For each $M \in \mathbb{N}$, there is a set of nodes $\{\theta_i^M\}_{i=1}^{M} \subset \Theta$ and an associated set of weights $\{\alpha_i^M\}_{i = 1}^M \subset \mathbb{R}$, such that for any continuous function $h:\Theta \to \mathbb{R}$,
\small
\begin{align*}
\int_{\Theta} h(\theta) d\theta  =  \lim_{M \to \infty} \sum_{i = 1}^M h(\theta^M_i) \alpha_i^M.
\end{align*}
\end{assumption}
\normalsize
This is the same as \cite{walton_IJC} Assumption 3.1; we include it for reference.

In additions to the assumptions of \cite{walton_IJC}, we also impose the following:

\begin{assumption} \label{regularityassumption}
The functions $f$ and $r$ are $C^1$.  The set $\Theta$ is compact and  $x_0:\Theta \mapsto \mathbb{R}^{n_x}$ is continuous.  Moreover, for the compact sets $X$ and $U$ defined in  \cite{walton_IJC}'s Assumptions 2.3 and 2.4, and for each $t \in [0,T]$, $\theta \in \Theta$, the Jacobians $r_x$ and $f_x$ are Lipschitz on the set $X \times U$, and the corresponding Lipschitz constants $L_r$ and $L_f$ are uniformly bounded in $\theta$ and $t$.  The function $F$ is $C^1$ on $X$ for all $\theta \in \Theta$; in addition, $F$ and $ F_x$ are continuous with respect to $\theta$. 
\end{assumption}

\subsection{Main Theorem Proof}

The theorem relies on the following  lemma:

\begin{lemma}\label{costate lemma} 
	Let $\{u_M\}$ be a sequence of optimal controls for Problem $\mathbf{P^M}$ with an accumulation point $u^\infty$ for the infinite set $V\subset \mathbb{N}$. Let $(x^\infty(t,\theta),\lambda^\infty(t,\theta))$ be the solution to the dynamical system:
\small
\begin{align}\label{costate accumulation points}
\begin{cases}
\dot{x}^\infty(t,\theta)=f(x^\infty(t,\theta),u^\infty(t),\theta) \\
\dot{\lambda}^\infty(t,\theta)= -\frac{\partial H(x^\infty(t,\theta),\lambda^\infty(t,\theta),u^\infty(t),t,\theta)}{\partial x}
\end{cases}
\end{align}
\vspace{-8pt}
\begin{align}
\begin{cases}
x^\infty(0,\theta) = x_0(\theta) \\
\lambda^\infty(T,\theta) = \frac{\partial F(x^\infty(T,\theta),\theta) }{\partial x}  
\end{cases}
\end{align}
\normalsize
where $H$ is defined as per Equation (\ref{first_Hamiltonian}), and let $\{(x_M(t,\theta),\lambda_M(t,\theta))\}$ for $M\in V$ be the sequence of solutions to the dynamical systems:
\small
\begin{align}\label{costate sequence}
\begin{cases}
\dot{x}_M(t,\theta)=f(x_M(t,\theta),u_M(t),\theta)\\
\dot{\lambda}_M(t,\theta)=  -\frac{\partial H(x_M(t,\theta),\lambda_M(t,\theta),u_M(t),t,\theta)}{\partial x}
\end{cases} 
\end{align}
\vspace{-8pt}
\begin{align}
\begin{cases}
x_M(0,\theta) = x_0(\theta) \\
\lambda_M(T,\theta) = \frac{\partial F(x_M(T,\theta),\theta) }{\partial x} 
\end{cases}
\end{align}
\normalsize

Then, the sequence $\{(x_M(t,\theta),\lambda_M(t,\theta))\}$ converges pointwise to  $(x^\infty(t,\theta),\lambda^\infty(t,\theta))$ and this convergence is uniform in $\theta$.
\end{lemma}
\subsubsection*{Proof}
\noindent The convergence of $\{x_M(t,\theta)\}$ is given by \cite{walton_IJC}, Lemmas 3.4, 3.5. The convergence of the sequence of solutions $\{\lambda_M(t,\theta)\}$ is guaranteed by the optimality of $\{u_M\}$. The convergence of  $\{\lambda_M(t,\theta)\}$ then follows the same arguments given the convergence of $\{x_M(t,\theta)\}$, utilizing the regularity assumptions placed on the derivatives of $F$, $r$, and $f$ with respect to $x$ to enable the use of Lipschitz conditions on the costate dynamics and transversality conditions. 
\begin{remark}\label{equivalence of sequence}
Note that $\lambda_M(t,\theta)$	is not a costate of Problem ${\mathbf{P^{ \lambda M}}}$, since it is a function of $\theta$. However, when $\theta = \theta_i^M$, then $\lambda_M(t,\theta_i^M) = \tilde{\lambda}_i^{M}(t)$, where $\tilde{\lambda}_i^{M}$ is the costate of Problem ${\mathbf{P^{ \lambda M}}}$ generated by the pair of solutions to Problem ${\mathbf{P^{M}}}$, $(\tilde{x}_i^{M}, u_M^\ast)$ . In other words, the function $\lambda_M(t,\theta)$ matches the costate values at all collocation nodes. Since these values satisfy the dynamics equations of Problem ${\mathbf{P^{ \lambda M}}}$, a further implication of this is that the values of $\lambda_M(t,\theta_i^M)$ produce feasible solutions to Problem ${\mathbf{P^{ \lambda M}}}$.
\end{remark}
\begin{remark}\label{equivalence of accumulation points}
Since the functions $\{(x_M(t,\theta),\lambda_M(t,\theta))\}$ obey the respective identities $x_M(t,\theta_i^M) = \tilde{x}_i^{M}(t)$ and $\lambda_M(t,\theta_i^M) = \tilde{\lambda}_i^{M}(t)$, their convergence to $(x^\infty(t,\theta),\lambda^\infty(t,\theta))$ also implies the convergence of the sequence of discretized primals and duals, $\{\tilde{X}_M\}$ and $\{\tilde{\Lambda}_M\}$, to accumulation points given by the relations 
\small
\[
\lim_{M\in V} \tilde{x}_i^{M}(t) = x^\infty(t,\theta_i^M), \quad
\lim_{M\in V} \tilde{\lambda}_i^{M}(t) = \lambda^\infty(t,\theta_i^M)
\]
\normalsize
\end{remark}	

We now prove Theorem \ref{hamiltonian convergence theorems}. Let $\{(x_M(t,\theta),\lambda_M(t,\theta))\}$ for $M\in V$ be the sequence of solutions defined by Equation \ref{costate sequence} and let $(x^\infty(t,\theta),\lambda^\infty(t,\theta))$ be the accumulation functions defined by Equation \ref{costate accumulation points}.  Incorporating Remarks \ref{equivalence of sequence} and \ref{equivalence of accumulation points}, we have:
\vspace{-5pt}
\small
\[
\lim_{M\in V} \tilde{H}^M(\tilde{X}_M, \tilde{\Lambda}_M, u_M,t)
=
\qquad \qquad \qquad \qquad \qquad \qquad \qquad \qquad \qquad 
\]
\vspace{-7pt}
\[
\lim_{M\in V}
\sum_{i=1}^M  \alpha_i^M \left[\tilde{\lambda}_i^M(t) f(\tilde{x}_i^M(t),u(t),\theta_i^M) +
\right.
\]
\vspace{-5pt}
\[
\left.
r(\tilde{x}_i^M(t), u(t),t,\theta_i^M) \right] 
\]
\vspace{-7pt}
\[
=
\lim_{M\in V}
\sum_{i=1}^M \alpha_i^M \left[ \lambda_M(t,\theta_i^M)f(x_M(t,\theta_i^M),u(t),\theta_i^M) + 
\right.
\]
\vspace{-5pt}
\[
\left. r(x_M(t,\theta_i^M), u(t),t,\theta_i^M) \right] 
\]
\normalsize
 
Due to the results of Lemma \ref{costate lemma}, and applying \cite{walton_IJC}'s Remark 1 on the convergence of the quadrature scheme for uniformly convergent sequences of continuous functions, we find that:
\vspace{-5pt}
\small
\[
\lim_{M\in V} \tilde{H}^M(\tilde{X}_M, \tilde{\Lambda}_M, u_M,t)
=
\qquad \qquad \qquad \qquad \qquad \qquad \qquad \qquad \qquad 
\]
\vspace{-5pt}
\[
\int_{\Theta} \left[ \lambda^\infty(t,\theta)f(x^\infty(t,\theta),u^\infty(t),\theta) + 
\right. \qquad \qquad \qquad \qquad \qquad \qquad  
\]
\vspace{-5pt}
\[
\left. r(x^\infty(t,\theta), u^\infty(t),t,\theta) \right]  d\theta 
= 
\mathbf{H}(x^\infty,\lambda^\infty,u^\infty,t)
\]
\normalsize
Thus proving the theorem.

\section*{Acknowledgment}

This work was supported in part by ONR SoA program and by NPS Cruser program. 
\ifCLASSOPTIONcaptionsoff
  \newpage
\fi

\bibliographystyle{IEEEtran}
\bibliography{references}  

\end{document}